# GROWTH OF THE BROWNIAN FOREST

By Jim Pitman[1] and Matthias Winkel[2]

*University of California, Berkeley and University of Oxford*

Trees in Brownian excursions have been studied since the late 1980s. Forests in excursions of Brownian motion above its past minimum are a natural extension of this notion. In this paper we study a forest-valued Markov process which describes the growth of the Brownian forest. The key result is a composition rule for binary Galton–Watson forests with i.i.d. exponential branch lengths. We give elementary proofs of this composition rule and explain how it is intimately linked with Williams' decomposition for Brownian motion with drift.

**1. Introduction.** Given $0 \leq \lambda < \mu$, the binary$(\lambda, \mu)$ forest is defined as a collection of independent binary Galton–Watson trees, with branching probability $(\mu - \lambda)/2\mu$, branch lengths independent exponential$(2\mu)$, planted into the positive real line at the points of a homogeneous Poisson process of rate $\mu - \lambda$.

We call the vertices on the positive real line *roots*. The unique branch connecting to a root is its *trunk*. We regard the forest as a plane forest: above any branch point we distinguish a left and a right *subbranch*. We call the positive real line the *forest floor*, and branch ends with no branching *leaves*. It is well known that the trees are critical for $\lambda = 0$ and subcritical for $\lambda > 0$, hence a.s. finite in either case.

Such random trees have been studied by a number of authors, including Neveu [21, 22, 23, 24], Le Gall [19], Geiger and Kersting [10], Shapiro [27] and Hobson [14]. These papers describe how in the critical case $\lambda = 0$ these trees may be found embedded in Brownian excursions, by a natural generalization of Harris's embedding of critical Galton–Watson trees with geometric offspring distribution in the excursions of a simple random walk

Received April 2004; revised November 2004.

[1]Supported in part by NSF Grants DMS-00-71448 and DMS-04-05779.

[2]Supported by Aon and the Institute of Actuaries.

*AMS 2000 subject classifications.* 60J65, 60J80.

*Key words and phrases.* Continuum random tree, forest-valued Markov process, binary Galton–Watson branching process, Brownian motion, Williams' decomposition.

[13]. We show in this paper how the random forests considered here can also be found embedded in Brownian paths in the subcritical case $\lambda > 0$, as increments of a forest-valued Markov process which describes the growth of the Brownian forest.

A common tool in studying trees and forests is searching trees by depth-first search (DFS) by which we mean

(a) starting at the root travelling continuously at unit speed along the branches;

(b) at each branch point first following the left subbranch, then the right subbranch, and then back the parent branch;

(c) at each leaf just following the branch back downward.

This generalizes to forests by starting to travel at location 0 on the forest floor in positive direction, searching all trees encountered, as described above.

Two binary forests $\mathcal{F}^{\lambda,\mu} \sim \mathrm{binary}(\lambda,\mu)$ and $\mathcal{F}^{\mu,\theta} \sim \mathrm{binary}(\mu,\theta)$ for $0 \leq \lambda < \mu < \theta$ can be composed as illustrated in Figure 1, by wrapping $\mathcal{F}^{\mu,\theta}$ (with thin lines) around the sides of $\mathcal{F}^{\lambda,\mu}$ (with thick lines).

DFS of $\mathcal{F}^{\lambda,\mu}$ at unit speed (up and down) maps time $[0,\infty)$ onto the graph of $\mathcal{F}^{\lambda,\mu}$. Identifying time $[0,\infty)$ with the original forest floor of $\mathcal{F}^{\mu,\theta}$, DFS maps the roots of $\mathcal{F}^{\mu,\theta}$ to points on the branches of $\mathcal{F}^{\lambda,\mu}$. We attach the corresponding trees of $\mathcal{F}^{\mu,\theta}$ at the images of their roots. Note that DFS visits every point on $\mathcal{F}^{\lambda,\mu}$ twice (except for the leaves and branch points), and we further specify that trees are attached to the left of a branch if visited for the first time, and to the right if visited for the second time. Every branch of length $t$ contributes a total length of $2t$ for the "planting" Poisson process. See [26] for a more formal construction of such trees and forests, in terms of which this composition operation can be rigorously defined.

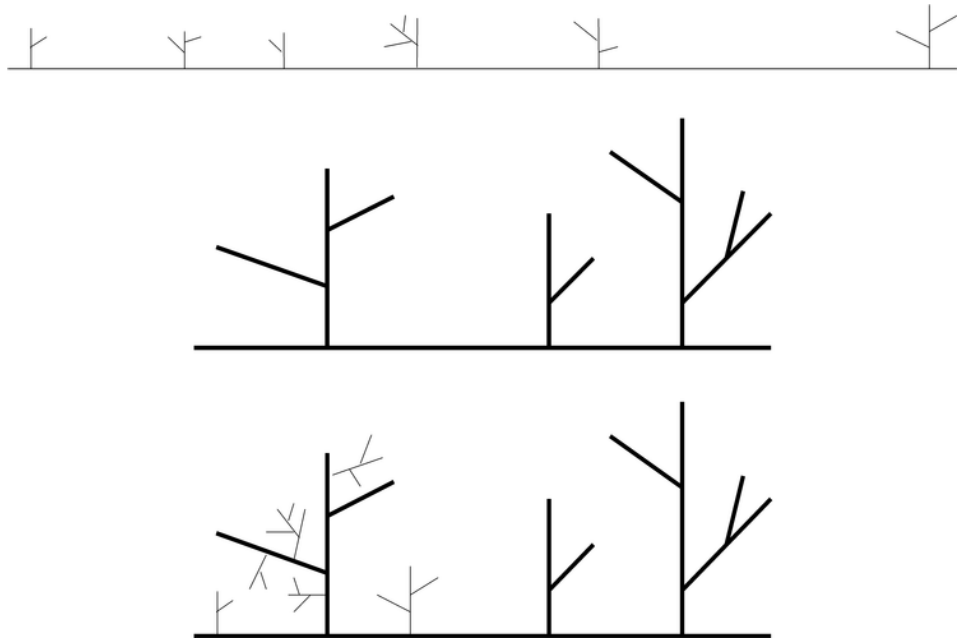

Fig. 1.  *Forest growth by wrapping one forest around another.*

THEOREM 1 (Forest composition rule). *For $0 \le \lambda < \mu < \theta$ the composition of a* binary$(\lambda,\mu)$ *forest with an independent* binary$(\mu,\theta)$ *forest is a* binary$(\lambda,\theta)$ *forest.*

In Section 2 we give a first elementary proof of this forest composition rule. Section 3 gives a different representation of binary$(\lambda,\mu)$ forests in terms of an alternating exponential random walk. Section 4 shows how this representation, first noted in the critical case $\lambda = 0$ by Le Gall [19] and Neveu and Pitman [23], establishes links with a Williams type decomposition for alternating exponential random walks. In Section 5 we show how the forest composition rule can be stated and proved in a Brownian forest setting, building on basic results about Brownian trees due to Neveu and Pitman [23, 24], Aldous [2] and Shapiro [27], and using links with Williams' [28] path decomposition of Brownian motion with drift.

For alternative approaches to Williams type decompositions we refer to [9, 11, 12, 16, 18, 20, 25]. For different tree-valued processes approaching Brownian trees we refer to [1, 3, 8, 22, 23].

The growth process described here is extended in [6] to forests with multiple branching, to allow growth of general Lévy forests formed from the Lévy trees studied by Duquesne and Le Gall [5].

## 2. Forest composition.

We refer once more to Figure 1 for a graphical representation of the forest composition rule, where the thick forest is meant to be the realization of a binary$(\lambda,\mu)$ forest, and the thin forest a binary$(\mu,\theta)$ forest. We will refer to the thick forest as the *black* forest, and think of the thin trees as *red* decoration. Clearly, each branch of the composite forest has a unique color. We can also associate colors with vertices. We need a convention for the color of vertices where red and black branches meet. The structure of the composed forest is such that these are all roots of the red forest, and we agree that these vertices are red.

We first reduce the study of composite forests to the study of the first tree in the composed forest.

LEMMA 1. *The trees of the composed forest are independent and identically distributed. They are separated by independent identically* exponential$(\theta - \lambda)$ *distributed spaces.*

PROOF. The first tree of the composed forest is clearly at an exponential($\theta - \lambda$) = min(exponential($\theta - \mu$), exponential($\mu - \lambda$)) position as the first of the first trees of the two forests.

Given that the trunk of the first tree is red, the residual black spacing is an exponential($\mu - \lambda$) variable independent of the first tree due to the lack of memory property. Hence the remaining black forest and clearly also the remaining red forest beyond the first tree are independent of the first tree and distributed like the whole forests.

Given that the trunk of the first tree is black, the residual last spacing of red trees on the black tree is an exponential($\theta - \mu$) variable independent of the first tree. Therefore the same conclusion holds.

An iteration concludes the proof. □

We next show that the branch lengths are independent of the shape, and i.i.d. exponential($2\theta$) distributed, where we mean by *shape* the combinatorial structure with unit branch lengths. In fact, we show the stronger statement

LEMMA 2. *The sequence of black branch lengths in the composite forest, enumerated by DFS, is an i.i.d. exponential($2\theta$) sequence, and so is the sequence of red branch lengths enumerated by DFS. The two sequences are independent and independent of the shapes of the trees.*

PROOF. We decompose the black trees by DFS and obtain a sequence $(A_m)_{m \geq 1}$ of i.i.d. exponential($2\mu$) branch lengths independent from the shapes of the black trees. We note that the red trees are independent of the Poisson process of their spacings. This Poisson process further decomposes each black branch into a geometric($\mu/\theta$) number such that each fragment $B_n$, $n \geq 1$, has an exponential($2\theta$) distribution.

The sequential construction shows that all fragment lengths $(B_n)_{n \geq 1}$ are independent and independent of the events $C_n$ that $B_n$ is associated with a fragment with a red tree at its far end, and the events $D_n$ that the red tree grows at the left (and not the right). We have $\mathbb{P}(C_n) = (\theta - \mu)/\theta$ and $\mathbb{P}(D_n | C_n) = 1/2$.

The statement on red branches and the remaining independence assertions follow from the corresponding properties of the red forest. □

COROLLARY 1. *The* binary($\lambda, \mu$) *forest decorated by an independent* binary($\kappa, \theta$) *forest has exponential branch lengths if and only if $\kappa = \mu$.*

PROOF. The same argument as for Lemma 2 shows that the black branches of the composite forest are i.i.d. exponential($2\theta - 2\kappa + 2\mu$) whereas the red branches are i.i.d. exponential($2\theta$). Therefore the length of the first trunk is a mixture with respective probabilities $(\mu - \lambda)/(\theta - \kappa + \mu - \lambda)$ and $(\theta - \kappa)/(\theta - \kappa + \mu - \lambda)$. This yields nonexponential densities unless $\kappa = \mu$. □

It remains to identify the shape of the first tree. Let us analyze the branching probabilities of the combined tree as decomposed above. Denote by $B_n^b$ the event that the $n$th black branch of the composite tree branches (enumeration by DFS). Then denote by $L_n^b$ and $R_n^b$ the colors of its left and right subbranches. Recall that we denote by $C_n^c$ the event that the $n$th *new* black branch is the uppermost fragment of an *old* black branch, and otherwise denote more specifically by $D_n \subset C_n$ the event that a red subtree is attached to the left of the black branch (as opposed to the right). *Old* and *new* refer to *before* and *after* the decoration by red trees.

COROLLARY 2. *The branching probabilities of red branches are given by*

| red/red | ∅ |
|---|---|
| $(\theta - \mu)/2\theta$ | $(\theta + \mu)/2\theta$ |

*The branching probabilities of (new) black branches are given by*

| red/black | black/red | black/black | ∅ |
|---|---|---|---|
| $(\theta - \mu)/2\theta$ | $(\theta - \mu)/2\theta$ | $(\mu - \lambda)/2\theta$ | $(\mu + \lambda)/2\theta$ |

*We interpret, for example, the third column as*

$$\mathbb{P}(B_n^b, L_n^b = \text{black}, R_n^b = \text{black}) = (\mu - \lambda)/2\theta.$$

*The composite trees show therefore a multitype branching structure with an absorbing red type and a black type that can produce both black and red children (but not two red children at the same time).*

PROOF. When enumerated by DFS, the branching events at the ends of the red branches are independent and identically distributed according to the red offspring distribution, as required.

Note that for any black fragment, enumerated by DFS, the events $C_n \cap D_n$, $C_n \cap D_n^c$, $C_n^c \cap B_n^b$ and $C_n^c \cap (B_n^b)^c$ describe the four possibilities for subbranches red/black, black/red, black/black and ∅, respectively. The argument of the proof of Lemma 2 yields the claimed probabilities and the independence of these events for different $n \geq 1$. □

Denote the color of the first trunk by $C_T$, its branching event by $B_T$ and the colors of its subbranches by $L_T$ and $R_T$.

Lemma 3. *The following conditional and unconditional probabilities hold*:

$$\mathbb{P}(C_T = \text{red}) = (\theta - \mu)/(\theta - \lambda),$$

$$\mathbb{P}(C_T = \text{black}) = (\mu - \lambda)/(\theta - \lambda),$$

$$\mathbb{P}(L_T = \text{red}|B_T) = \mathbb{P}(R_T = \text{red}|B_T) = (\theta - \mu)/(\theta - \lambda),$$

$$\mathbb{P}(L_T = \text{red}, R_T = \text{red}|B_T) = (\theta - \mu)^2/(\theta - \lambda)^2,$$

$$\mathbb{P}(C_T = \text{red}|B_T^c) = (\theta^2 - \mu^2)/(\theta^2 - \lambda^2).$$

Proof. The first two probabilities are deduced from the laws of the spacing variables.

For the others we calculate

$$\mathbb{P}(L_T = \text{red}|B_T) = \mathbb{P}(B_T)^{-1}\mathbb{P}(B_T, L_T = \text{red})$$
$$= \frac{2\theta}{\theta - \lambda}\left(\frac{\theta - \mu}{\theta - \lambda}\frac{\theta - \mu}{2\theta} + \frac{\mu - \lambda}{\theta - \lambda}\frac{\theta - \mu}{2\theta}\right) = \frac{\theta - \mu}{\theta - \lambda}$$

by splitting the event into red and black trunk, and finally

$$\mathbb{P}(L_T = \text{red}, R_T = \text{red}|B_T) = \frac{2\theta}{\theta - \lambda}\ \frac{\theta - \mu}{\theta - \lambda}\frac{\theta - \mu}{2\theta} = \left(\frac{\theta - \mu}{\theta - \lambda}\right)^2,$$

$$\mathbb{P}(C_T = \text{red}|B_T^c) = \mathbb{P}(B_T^c)^{-1}\mathbb{P}(B_T^c|C_T = \text{red})\mathbb{P}(C_T = \text{red})$$
$$= \frac{2\theta}{\theta + \lambda}\frac{\theta + \mu}{2\theta}\ \frac{\theta - \mu}{\theta - \lambda}.\qquad\square$$

Corollary 3. *Given the first tree has two first subbranches, the colors $C_L$ and $C_R$ of the two subbranches are independent and have the same distribution as the color $C_T$ of the trunk.*

*Given the first tree has no branches apart from the trunk, the distribution of the color of the trunk is given by $p_r = (\theta^2 - \mu^2)/(\theta^2 - \lambda^2)$ and $p_b = (\mu^2 - \lambda^2)/(\theta^2 - \lambda^2)$.*

This suggests that the colors of the leaves are i.i.d. according to these $p_r$ and $p_b$. To prove Theorem 1 we establish the stronger fact

Proposition 1. *The composite forest has the same law as a* binary$(\lambda, \theta)$ *forest whose leaves are colored by i.i.d. variables with law given by $p_r$ and $p_b$, with further coloring by the following rules. Start with a red tree and color black the subtree spanned by the root and the black leaves* (*if any*).

Proof. By Lemmas 1 and 2 it is sufficient to check the law of the shape and coloring of the first tree of the composite forest. Call $\mathcal{T}_1$ this tree with branch lengths replaced by unit lengths.

Call $\tilde{\mathcal{T}}_1$ the first tree of a binary$(\lambda, \theta)$ forest, colored according to the rules stated in the proposition, with unit branch lengths.

For any colored tree $t$ with $k \geq 0$ black leaves and $n - k$ red leaves, it is an elementary combinatorial fact that there are $k$ inner vertices (including the root unless $k = 0$) that have only black subbranches and $n - k$ inner vertices that have at least one red subbranch.

Fix a number $n \geq 1$ of leaves, the shape of the tree, a number $k$ of black leaves, the positions of the black leaves. We show that the probabilities that $\mathcal{T}_1$ and $\tilde{\mathcal{T}}_1$ are this tree $t$, are the same.

For $k = 0$ we calculate from the probabilities established in Corollary 2

$$\mathbb{P}(\mathcal{T}_1 = t) = \frac{\theta - \mu}{\theta - \lambda} \left( \frac{\theta - \mu}{2\theta} \right)^{n-1} \left( \frac{\theta + \mu}{2\theta} \right)^n,$$

$$\mathbb{P}(\tilde{\mathcal{T}}_1 = t) = \left( \frac{\theta - \lambda}{2\theta} \right)^{n-1} \left( \frac{\theta + \lambda}{2\theta} \right)^n \left( \frac{(\theta - \mu)(\theta + \mu)}{(\theta - \lambda)(\theta + \lambda)} \right)^n,$$

where the first product consists of the probabilities to obtain a red first trunk, $n - 1$ branchings and $n$ leaves (no branching) in the right order which is determined by the given shape. The second product consists of $n - 1$ branchings and $n$ leaves for the shape and $n$ red leaves for the coloring. The two products are equal by obvious cancellation.

Similarly for $k > 0$

$$\mathbb{P}(\mathcal{T}_1 = t) = \frac{\mu - \lambda}{\theta - \lambda} \left( \frac{\theta - \mu}{2\theta} \right)^{n-k} \left( \frac{\mu - \lambda}{2\theta} \right)^{k-1} \left( \frac{\theta + \mu}{2\theta} \right)^{n-k} \left( \frac{\mu + \lambda}{2\theta} \right)^k,$$

$$\mathbb{P}(\tilde{\mathcal{T}}_1 = t) = \left( \frac{\theta - \lambda}{2\theta} \right)^{n-1} \left( \frac{\theta + \lambda}{2\theta} \right)^n \left( \frac{(\theta + \mu)(\theta - \mu)}{(\theta + \lambda)(\theta - \lambda)} \right)^{n-k} \left( \frac{(\mu + \lambda)(\mu - \lambda)}{(\theta + \lambda)(\theta - \lambda)} \right)^k,$$

where the first product contains the probabilities of a black first trunk, $n - k$ branchings leading to at least one red subbranch, $k - 1$ branchings leaving no red subbranch, $n - k$ red leaves (no branching) and $k$ black leaves. Note that it is important that the $n - k$ branchings leading to at least one red subbranch all branch in the required way with the same probability $(\theta - \mu)/2\theta$ regardless of whether the parent branch is red or black and in the latter case whether the red subbranch is to the left or to the right. The second product is the same as above with the obvious changes in the coloring probabilities. Again, the two products are equal by cancellation. $\quad\square$

**3. Equivalent construction of binary$(\lambda, \mu)$ forests.** Let $(F_n)_{n \geq 0}$ and $(R_n)_{n \geq 1}$ be two independent sequences of i.i.d. exponential random variables with parameters $\mu - \lambda$ and $\mu + \lambda$, respectively. We construct a random forest by the following rules, that are best explained looking at a picture, for example, the first forest in Figure 1. A branch of length $R_1$ rises from location $F_0$ on the forest floor to what is to become the leftmost leaf of the first tree. At length $F_1$ below this leaf, another branch commences. If $F_1 > R_1$, this is to be understood as beginning a new tree on the forest floor at distance $F_1 - R_1$ to the right of the (root of the) first tree. For $n \geq 2$, in either situation, at the respective location, a branch of length $R_n$ rises to the $n$th leaf. Following the current tree downward (and along the forest floor, if necessary) from this leaf until $F_n$ below, we reach the location where the next branch (or tree) is to be grown.

$(F_n)_{n \geq 0}$ and $(R_n)_{n \geq 1}$ are the falls and rises in the so-called Harris path associated with the forest which is constructed as the distance from the root when travelling at unit speed along the branches of the tree, at each branching inductively following first the left subbranch, then the right subbranch, then back downward; at each leaf just going back downward (depth first search). This defines the Harris path as a process in continuous time with slopes $\pm 1$. When restricting to the discrete time set of changes in direction, no information is lost, and this process is called the Harris walk. It is a random walk $(H_n)_{n \geq 0}$ with $H_0 = 0$, $H_1 = -F_0$ and independent alternating increments

$$H_{2n} - H_{2n-1} = R_n, \qquad H_{2n+1} - H_{2n} = -F_n, \qquad n \geq 1,$$

where $R_n \sim$ exponential$(\mu + \lambda)$ and $F_n \sim$ exponential$(\mu - \lambda)$.

PROPOSITION 2. *The forest constructed from the alternating walk* $(H_n)_{n \geq 0}$ *with independent* exponential$(\mu - \lambda)$ *falls and* exponential$(\mu + \lambda)$ *rises is a* binary$(\lambda, \mu)$ *forest.*

PROOF. Le Gall [19] proves this for $\lambda = 0$ by splitting the associated Harris walk at its infimum. See also [10]. The proof for $\lambda = 0$ can be adapted to the case $\lambda > 0$, or that case can be deduced by density calculations between the laws of the first trees for $\lambda > 0$ and for $\lambda = 0$. We present here a direct method that only involves calculations with exponential variables.

The key fact is that for $R_1 \sim$ exponential$(\mu + \lambda)$ and $F_1 \sim$ exponential$(\mu - \lambda)$ independent, $L_1 = \min(R_1, F_1) \sim$ exponential$(2\mu)$ gives the required branch length. By careful use of the lack of memory property of the exponential distribution, $L_1$ is independent from $(M_1, O_1)$ where $M_1 = \mathbb{1}_{\{L_1 = F_1\}}$ and $O_1 = \max(R_1, F_1) - L_1$. $\{M_1 = 1\}$ is a branching event with $\mathbb{P}(M_1 = 1) = (\mu - \lambda)/2\mu$, and given $M_1 = 1$, $O_1 \sim$ exponential$(\mu + \lambda)$ has the rise distribution; given $M_1 = 0$, $O_1 \sim$ exponential$(\mu - \lambda)$ has the fall distribution. The overshoot $O_1$ will further contribute either as spacing between trees or for further comparisons with other exponential variables.

More precisely, we proceed inductively as follows. We set up an empty LIFO stack (last in first out) that will store overshoots with rise distribution. Queueing techniques have been popular in branching processes ever since the early work of Kendall [15]. The above first step defines two independent random variables $L_1$ and $M_1$, and an overshoot $O_1$ whose conditional distribution depends on $M_1$. We also record $N_1 = 1$ to keep track of how many original rises $R_n$ and falls $F_n$ we have used. Stack size is $S_1 = 0$.

Our induction hypothesis is now that we have constructed two independent sequences of i.i.d. random variables $(L_m)_{1 \le m \le k}$ and $(M_m)_{1 \le m \le k}$ and an overshoot $O_k$ whose conditional distribution depends only on $M_k$. Furthermore, assume that given $S_k$, the stack is independent of $(L_m)_{1 \le m \le k}$, $(M_m)_{1 \le m \le k}$ and $O_k$ and contains $S_k$ i.i.d. variables $\sim$ exponential$(\mu + \lambda)$.

If $M_k = 1$, put $O_k$ on the stack, and construct $(L_{k+1}, M_{k+1}, O_{k+1})$ from the next rise $R_{N_{k+1}}$ and fall $F_{N_{k+1}}$ where $N_{k+1} = N_k + 1$, $S_{k+1} = S_k + 1$.

If $M_k = 0$ and $S_k > 0$, take $\tilde{R}$ from the stack and construct $(L_{k+1}, M_{k+1}, O_{k+1})$ from $\tilde{R}$ and $\tilde{F} = O_k$. Note that by hypothesis, given $M_k = 0$ and $S_k > 0$, $\tilde{R}$ is independent of $\tilde{F}$. Put $N_{k+1} = N_k$, $S_{k+1} = S_k - 1$.

If $M_k = 0$ and $S_k = 0$, a tree has been finished, and $O_k$ is the spacing to the next tree. Construct $(L_{k+1}, M_{k+1}, O_{k+1})$ from $R_{N_{k+1}}$ and $F_{N_{k+1}}$ where $N_{k+1} = N_k + 1$, $S_{k+1} = S_k = 0$.

Then, $(L_{k+1}, M_{k+1}, O_{k+1})$ is conditionally independent from $(L_m)_{1 \le m \le k}$ and $(M_m)_{1 \le m \le k}$ given $M_k$ and $S_k$, and its conditional distribution does not depend on $M_k$ and $S_k$, so they are unconditionally independent. Therefore $(L_m)_{1 \le m \le k+1}$ and $(M_m)_{1 \le m \le k+1}$ are independent sequences of i.i.d. random variables and the conditional distribution of $O_{k+1}$ only depends on $M_{k+1}$. By hypothesis, the stack was conditionally independent of $(L_m)_{1 \le m \le k}$, $(M_m)_{1 \le m \le k}$ and $O_k$ given $S_k$ and consisted of i.i.d. rises. Therefore, now, given $M_k$ and $S_k$ (hence also given $S_{k+1}$), the stack is conditionally independent of $(L_m)_{1 \le m \le k}$, $(M_m)_{1 \le m \le k}$, $(L_{k+1}, M_{k+1}, O_{k+1})$. Since the conditional distribution of the stack and $(L_m)_{1 \le m \le k+1}$, $(M_m)_{1 \le m \le k+1}$ and $O_{k+1}$ only depends on $S_{k+1}$, independence holds given $S_{k+1}$. By hypothesis, the conditional distribution of the stack given $S_k$ (and the conditionally independent $M_k$) was that of $S_k$ i.i.d. variables $\sim$ exponential$(\mu + \lambda)$. Therefore, in each of the three cases, the conditional distribution of the stack given $M_k$ and $S_k$ (hence also given $S_{k+1}$) is that of $S_{k+1}$ i.i.d. variables $\sim$ exponential$(\mu + \lambda)$. Specifically, note in the first case that the new variable $O_k$ put onto the stack, jointly with $M_k$ is conditionally independent of the stack given $S_k$, and given also $M_k = 1$, has an exponential$(\mu + \lambda)$ distribution, as required. This completes the induction step.

The induction shows that the procedure yields independent sequences $(L_m)_{m \ge 1}$ and $(M_m)_{m \ge 1}$ with $L_m \sim$ exponential$(2\mu)$ and $\mathbb{P}(M_m = 1) = 1 - \mathbb{P}(M_m = 0) = (\mu - \lambda)/2\mu$. In the construction described at the beginning of this section, $L_m$, $m \ge 1$, are an enumeration of the branch lengths whereas $M_m$, $m \ge 1$, determine the shapes of the binary trees.

The distribution of the numbers of leaves of a tree under this construction is easily identified with the one of a binary$(\lambda, \mu)$ tree since, in both cases, a tree is finished once there have been more no-branchings than branchings (here more $M_m = 1$ than $M_m = 0$). Also clearly, in both cases, the shape of a tree with $n$ leaves is uniformly distributed (coded by the order of branchings and no-branchings). Just note that, combinatorially, the tree shapes are different for each two different (admissible) orders of $M_1, \ldots, M_{n-1}$ because the associated Harris walks (with unit rises and falls) are then different and code different tree shapes.    □

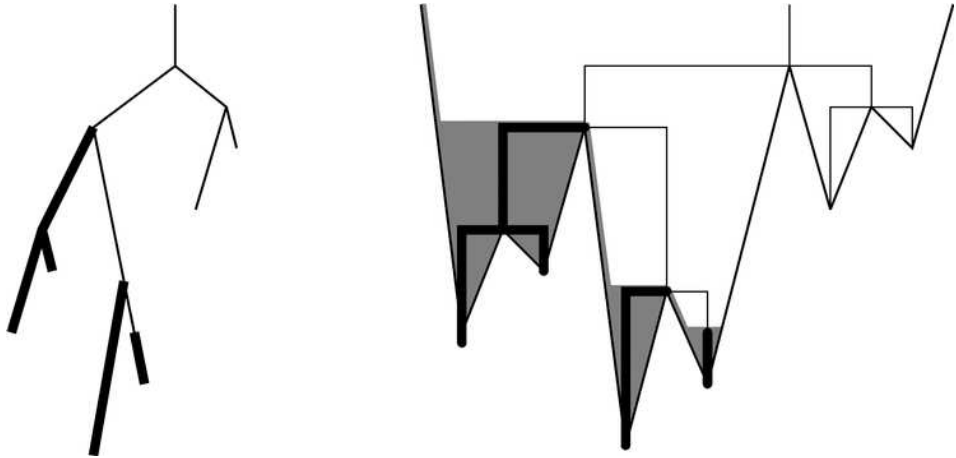

Fig. 2.    *The tree in a Harris walk built from the leaves.*

The proof suggests building the trees from the leaves. From left to right, each leaf is placed and gradually connected to subtrees all of whose leaves have been placed. In the associated Harris walk $(H_n)_{n \geq 0}$, leaves are the local maxima. The order in which the tree is constructed can be illustrated in the associated Harris path with the tree drawn underneath, or even better in $(-H_n)_{n \geq 0}$. Imagine the path $(-H_n)_{n \geq 0}$ as a (two-dimensional) container in which we pour water from the origin, see Figure 2. The rest of the description is what naturally happens. The first part to fill corresponds to the first leaf. When this is filled, water pours into the second leaf, and so on. If these two leaves form a subtree, their common parent branch is the next to fill. If not, the next leaf comes next and will be connected later to the previous leaves in reverse order.

In Knuth's terminology, our (binary) tree is drawn in "Postorder" when traversing it by DFS since a node is only recognized as such when it is visited for the third (and last) time, as opposed to the usual "Preorder" (if nodes are recognized while rising, as in the Galton–Watson tree construction in the Introduction) or "Inorder" (if nodes are not known while rising but recognized when rising into the second subbranch, as when drawing the tree directly from rises and falls); see [17], Section 2.3.1.

**4. Williams type decomposition of Harris walks.** We study here the Harris walks of the three forests in the forest composition. As we have seen in the previous section, each of them is essentially an alternating exponential random walk with the mean rise less than or equal to the mean fall, so that their lim inf is $-\infty$. In the sequel we shall switch between Harris walks (discrete time) and Harris paths (continuous time, slopes $\pm 1$) since the walks are nicer stochastic processes to formulate a Williams type decomposition and easier when passing to limits in the next section, whereas the paths are more suitable to describe the forest composition and technically derive Williams' decomposition.

Three subsections treat preliminaries, derivations of the Williams type decomposition for Harris walks from the forest composition rule, and vice versa, respectively.

4.1. *Preliminaries on forest structure and Harris paths.* The first (deterministic) result on Harris paths of composite forests does not depend on the special (stochastic) structure of binary$(\lambda, \mu)$ forests. Its proof and also extensions to forests with nonbinary (but finite) trees and locations of multiple trees are straightforward.

PROPOSITION 3. *The Harris path $H^c$ of the composite forest is equal to the Harris path $H^b$ of the black forest with excursions inserted at the times of the arrival process of red trees where each excursion is the Harris path of the corresponding tree (without the extension into the negative half line that represents the spacing to the next tree).*

*More precisely, denote by $(T_m)_{m \geq 1}$ the locations of red trees on their forest floor, $t \to N_t$ the associated counting process. Let $s \mapsto H^{(m)}(s)$ be the Harris path of the $m$th red tree with lifetime $L_m = \inf\{s > 0 \colon H^{(m)}(s) = 0\}$. Then we have*

$$H^c\left(t + \sum_{m=1}^{N_t} L_m\right) = H^b(t),$$

$$H^c\left(T_n + \sum_{m=1}^{n-1} L_m + s\right) = H^b(T_n) + H^{(n)}(s),$$

*for all $t \geq 0$, $n \geq 1$ and $0 \leq s \leq L_n$.*

The second result studies further the red excursions inserted into the black Harris path.

LEMMA 4. *For a binary$(\mu, \theta)$ tree, $0 \leq \mu < \theta$, the Harris path $H(s)$ with lifetime $L = \inf\{s > 0 \colon H(s) = 0\}$ is reversible, that is,*

$$(H(s))_{0 \leq s \leq L} \sim (H(L - s))_{0 \leq s \leq L}.$$

PROOF.   This is obvious if we use the original representation of binary$(\mu, \theta)$ trees in terms of lengths and branchings (Proposition 2), which is symmetric under exchanging left and right.   □

The two results of this subsection are important when studying the relationship between Williams type decompositions of Harris walks and forest composition. We will need them, both to derive the decomposition from the forest composition and vice versa.

4.2. *Williams' decomposition from the forest composition rule.*

THEOREM 2 (Williams' decomposition for Harris walks).   *Let* $(H_n)_{n \geq 0}$ *be an alternating walk with independent* exponential$(\theta - \lambda)$ *falls and* exponential$(\theta + \lambda)$ *rises, for some* $0 \leq \lambda < \theta$, *and* $N \sim$ geometric$(1 - q)$ *an independent geometric random variable. Define the (absolute) height* $F = -\min\{H_n : 0 \leq n \leq 2N\}$ *and time* $M$, $F = -H_M$, *of the minimum before* $2N$, *and the terminal height* $R = H_{2N} + F$ *above the minimum.*

*Then the walk* $(H_n)_{0 \leq n \leq 2N}$ *decomposes into two independent parts* $(H_n)_{0 \leq n \leq M}$ *and* $(H_{M+n} - H_M)_{0 \leq n \leq 2N - M}$. *Furthermore, if we put* $\mu^2 = q\lambda^2 + (1 - q)\theta^2$,

$$F \sim \text{exponential}(\mu - \lambda),$$

$$R \sim \text{exponential}(\mu + \lambda),$$

$$(H_n)_{0 \leq n \leq M} \sim (\tilde{H}_n \vee (-F))_{0 \leq n \leq \tau(-F)},$$

$$(H_{2N} - H_{2N-n})_{0 \leq n \leq 2N - M} \sim (\tilde{H}_n \vee (-R))_{0 \leq n \leq \tau(-R)},$$

*where* $\tilde{H}$ *is independent of* $F$ *and* $R$ *and distributed as* $H$ *but with* exponential$(\theta - \mu)$ *falls and* exponential$(\theta + \mu)$ *rises, and* $\tau(-x) = \inf\{n \geq 0 : \tilde{H}_n < -x\}$.

PROOF.   By Proposition 2, the exponential falls and rises define a binary$(\lambda, \theta)$ forest. Furthermore, this is the setting of the forest composition where $N$ is the number of leaves before and including the first black leaf. This can be seen from Proposition 1 since the leaf colors in the combined forest are i.i.d., red with probability $q = p_r$, so that the first black leaf occurs at a geometric$(1 - q)$ position independent of the binary$(\lambda, \theta)$ forest.

We can identify the other quantities: $F$ and $R$ are the first fall and rise in the Harris walk of the black forest since by Proposition 3 the red trees insert only (positive) excursions into the black Harris walk, so that the minimum corresponds to a black change from fall to rise, by definition of $N$ the first and only such change before $2N$ (recall that the walk does two steps for every leaf). $M$ is the position of this minimum in the combined Harris walk.

The independence statement is obvious from the definition of forest composition and Proposition 3, since independent geometric numbers of independent excursions from red trees are inserted into the first fall and first rise of the black Harris path. The law of the first part is clear, too, since before the first rise of the black walk, the combined forest coincides with the red forest.

The law of the second part is obtained by time reversal. In fact, this part of the Harris walk is not simple when considered in the original sense, since the red excursions that start with a rise and end with a fall are inserted into a rise, so that the combined walk has nonexponential laws for the total rise before and at the beginning of each excursion, and again a nonexponential fall at the end of the excursion because it is an undershoot of the passage event into the negative half line. If we reverse time, however, the Poisson locations of trees are still Poisson locations, and the excursions still have the same laws, by the preceding lemma. Furthermore, the excursions are now embedded into a fall, and the same argument as for the first part identifies the law. $\square$

4.3. *The forest composition rule from Williams' decomposition.* Alternatively, Theorem 2 can be proved using Theorem 2 in Kersting and Memisoglu [16]. Their theorem is a general splitting result for Markov chains, here $X_n = (H_{2n}, H_{2n+1})$ killed at $N \sim \text{geometric}(1-q)$, at times

$$T = \sup\{n \geq 0 : h(X_m) < h(X_n) \leq Y \text{ for all } m < n\},$$

where $h$ is nonnegative and harmonic for $X$; here $h(u,v) = e^{-(\mu-\lambda)v}$ and $Y^{-1} \sim \text{uniform}(0,(h(0))^{-1})$ independent. With this choice of $h$, we have $h(0) = 1$ and identify

$$T = \sup\{n \geq 0 : H_{2m+1} > H_{2n+1} \geq -F \text{ for all } m < n\} = M,$$

$$F = \frac{1}{\mu - \lambda} \log(Y) \sim \text{exponential}(\mu - \lambda).$$

Kersting and Memisoglu give the distribution of the first subpath in terms of the $h$-transformed semigroup of $X$, and in our special case, this is easily seen to correspond to alternating walks with exponential$(\theta - \mu)$ falls and exponential$(\theta + \mu)$ rises. To also identify the distribution of the second subpath, a further time-reversal argument is needed, like in our proof.

Actually, the independence statement in Theorem 2 does not depend on the step distribution. It follows from the corresponding result for the stopped random walk $S_n = H_{2n}$, $0 \leq n \leq N$, in the augmented filtration $(\mathcal{G}_n)_{n \geq 0}$ of $(S_n)_{0 \leq n \leq N}$ that contains $H_{2n-1} \in \mathcal{G}_n$, see [20]. This argument can be expanded to prove Theorem 2 without using the forest composition rule.

In this special situation, another closely related approach is possible, since the renewal process of descending ladder heights has exponential increments, and excursions of the (discrete-time) process $H_n$ above the infimum are naturally described by a Poisson point process indexed by the negative spatial axis. General splitting of Poisson point processes as in this setting has been discussed by Greenwood and Pitman [12]. The following will lead to a stronger version of Theorem 2 that highlights the excursions above an alternating lower envelope, which we can then use to give a second proof of the forest composition rule, Theorem 1.

LEMMA 5. *Let* $(H_n)_{n \geq 0}$ *be an alternating walk with independent* exponential$(\theta - \lambda)$ *falls and* exponential$(\theta + \lambda)$ *rises. Define* $T_1 = 1$, $V_m = -H_{T_m}$, $T_{m+1} = \inf\{n \geq 0 : H_n < -V_m\}$, $m \geq 1$. *Then* $(V_m)_{m \geq 0}$ *is a renewal process with* exponential$(\theta - \lambda)$ *step distribution,* $C_t = \sum_{m \geq 1} \mathbb{1}_{\{V_m \leq t\}}$ *is a Poisson process with intensity* $\theta - \lambda$ *and*

$$\varepsilon_{V_m} = (H_{T_m + n} - H_{T_m})_{0 \leq n < T_{m+1} - T_m}, \qquad \varepsilon_t = \partial \qquad otherwise,$$

*defines a Poisson point process* $(\varepsilon_t)_{t \geq 0}$ *on the path space* $\bigcup_{n \geq 0}[0, \infty)^n$, $[0, \infty)^0 := \{\partial\}$ *is used as cemetery, with intensity measure*

$$\nu^{\lambda, \theta} = (\theta - \lambda)\mathbb{P}((H_{1+n} - H_1)_{0 \leq n < T_2 - 1} \in \cdot).$$

PROOF. Clearly $-H_1 \sim$ exponential$(\theta - \lambda)$. The overshoot of $H$ into $(-\infty, H_1)$ has an exponential$(\theta - \lambda)$ distribution, for example, because on $\{T_2 = 2n + 1\}$, $H_{2n+1} - H_1$ has this distribution, for all $n \geq 1$, by the lack of memory property of the exponential distribution. Also, the overshoot is easily seen to be independent of $(H_n)_{1 \leq n < T_2}$.

The remaining assertions are analogous to the theory of ladder events for random walks and their excursions above the minimum. The independence of ladder times and ladder heights and the exponential step distributions of ladder heights here allow us to express this in terms of Poisson point processes. □

This lemma establishes the bridge to the splitting result of [12] that we use to establish the Williams decomposition for Harris walks without using the forest composition rule. These methods are more sophisticated than our elementary approach via the forest composition, yet less sophisticated than an alternative derivation from embedding into Brownian paths by Poisson sampling, see the proof of Lemma 9. More precisely, we establish the following more general result.

PROPOSITION 4. *Let* $(H_n)_{n \geq 0}$ *be an alternating walk with independent* exponential$(\theta - \lambda)$ *falls and* exponential$(\theta + \lambda)$ *rises, for some* $0 \leq \lambda < \theta$, and $N_0 = 0$, $N_k - N_{k-1} \sim \text{geometric}(1-q)$, $k \geq 1$, independent time lags. Denote by

$$M_k = \operatorname*{arg\,min}_{2N_k \leq n \leq 2N_{k+1}} H_n, \qquad F_k = H_{2N_k} - H_{M_k},$$

$$R_{k+1} = H_{2N_{k+1}} - H_{M_k}, \qquad k \geq 0,$$

the intertwining minima and the associated falls and rises. Then the processes $\acute{\varepsilon}^{(0)} := (\varepsilon_t)_{0 \leq t < F_0}$ from Lemma 5, and $(\acute{\varepsilon}_t^{(k)})_{0 \leq t < R_k}$ and $(\grave{\varepsilon}_t^{(k)})_{0 \leq t < F_k}$, defined analogously from

$$\acute{H}_n^{(k)} = H_{2N_k - n} - H_{2N_k}, \qquad 0 \leq n \leq 2N_k - M_{k-1},$$

$$\grave{H}_n^{(k)} = H_{2N_k + n} - H_{2N_k}, \qquad 0 \leq n \leq M_k - 2N_k,$$

$k \geq 1$, are independent Poisson point processes with intensity measure $\nu^{\mu,\theta}$ killed at independent times $F_k \sim \text{exponential}(\mu - \lambda)$, $k \geq 0$, and $R_k \sim \text{exponential}(\mu + \lambda)$, $k \geq 1$.

PROOF (AND SECOND PROOF OF THEOREM 2). Mark each rise of the Harris walk $(H_n)_{n \geq 0}$ with probability $1 - q$ according to a sequence of independent Bernoulli variables $B_{2k}$ with

$$\mathbb{P}(B_{2k} = 0) = q, \qquad \mathbb{P}(B_{2k} = 1) = 1 - q, \qquad k \geq 1.$$

The first mark will be at $N = \inf\{k \geq 1 : B_{2k} = 1\} \sim \text{geometric}(1-q)$.

Following the lines of the proof of Lemma 2.1 of [12], we take the excursion process $\varepsilon$ studied in the previous lemma and record as a mark $\Gamma_{V_m} = \partial$ if $\varepsilon_{V_m}$ does not contain a mark, that is, if $B_{T_m + 2k-1} = 0$ for all $2k - 1 < L_m = \inf\{n \geq 0 : H_{T_m+n} - H_{T_m} < 0\} = T_{m+1} - T_m$. Otherwise put

$$\Gamma_{V_m} = \inf\{1 \leq 2k - 1 < L_m : B_{T_m + 2k-1} = 1\},$$

so that $2\Gamma_{V_m} - 1$ is the location of the first mark in excursion $\varepsilon_{V_m}$, if any. Put $\Gamma_t = \partial$ otherwise.

Then $(\varepsilon, \Gamma)$ is a Poisson point process and the first marked point is at $F = \inf\{t \geq 0 : \Gamma_t \neq \partial\}$. If $(G_m)_{m \geq 1}$ is a sequence of independent geometric$(1-q)$ variables, note that $(\varepsilon, \Gamma)$ has the same distribution as $(\varepsilon, \Gamma^*)$ where $\Gamma_{V_m}^* = G_m$ if $2G_m - 1 < L_m$, $\Gamma_t^* = \partial$ otherwise. By Lemma 3.3 of [12]:

(a) The law of $F$ is exponential with rate

$$r = \sum_{1 \leq n < \infty} (1 - q^n) \nu_n([0, \infty)^{2n}),$$

where $\nu_n = \nu(\cdot \cap [0, \infty)^{2n})$ is the characteristic measure of excursions of length $2n$.

(b) $(\varepsilon_t)_{0 \le t < F} \sim (\tilde{\varepsilon}_t)_{0 \le t < F}$ where $\tilde{\varepsilon}$ is an independent Poisson point process with characteristic measure

$$\sum_{1 \le n < \infty} q^n \nu_n$$

since an excursion of length $2n$ is not marked with probability $q^n$.

(c) The law of $(\varepsilon_F, \Gamma_F)$ is

$$\frac{1}{r} \sum_{1 \le n < \infty} (1 - q^n) \nu_n \otimes \gamma_n,$$

where $\gamma_n(\{j\}) = q^{j-1}(1-q)/(1-q^n)$ since an excursion of length $2n$ is marked with probability $1 - q^n$, and given that it is marked, the position of the first mark has a truncated geometric distribution $\gamma_n$.

(d) $F$ and $(\varepsilon_t)_{0 \le t < F}$ are independent of $(\varepsilon_F, \Gamma_F)$.

Therefore, we only need to compute a few quantities. We will see now that, as stated as part of Theorem 2, $\mu^2 = q\lambda^2 + (1-q)\theta^2$, that is, $q = (\theta^2 - \mu^2)/(\theta^2 - \lambda^2)$, is the $\mu$-value corresponding to $q$ that leads to the correct parameter of $F$. We use the walk-tree correspondence (Proposition 2, which we proved independently of the forest composition rule) to calculate

$$r = (\theta - \lambda) \sum_{1 \le n < \infty} \left(1 - \left(\frac{\theta^2 - \mu^2}{\theta^2 - \lambda^2}\right)^n\right) C_n \left(\frac{\theta - \lambda}{2\theta}\right)^{n-1} \left(\frac{\theta + \lambda}{2\theta}\right)^n$$

$$= \theta - \lambda - (\theta - \mu) \sum_{1 \le n < \infty} C_n \left(\frac{\theta - \mu}{2\theta}\right)^{n-1} \left(\frac{\theta + \mu}{2\theta}\right)^n = \mu - \lambda,$$

where $C_n$ is the number of shapes of binary trees with $n$ leaves, and $p^{n-1}(1-p)^n$ is the probability that each of these occurs where $p$ is the respective branching probability in a binary$(\lambda, \theta)$ forest and a binary$(\mu, \theta)$ forest.

For binary Galton–Watson trees, the conditional shape distribution is uniform on the $C_n$ possibilities given the number $n$ of leaves. Since shape and branch lengths are independent, by definition, we deduce that $\nu_n/\nu_n([0, \infty)^{2n})$ does not depend on $\lambda$; therefore, a variant of the previous computation shows that

$$\sum_{1 \le n < \infty} q^n \nu_n = \sum_{1 \le n < \infty} \left(\frac{\theta^2 - \mu^2}{\theta^2 - \lambda^2}\right)^n (\theta - \lambda) C_n \left(\frac{\theta - \lambda}{2\theta}\right)^{n-1} \left(\frac{\theta + \lambda}{2\theta}\right)^n \frac{\nu_n}{\nu_n([0, \infty)^{2n})}$$

$$= (\theta - \mu) \mathbb{P}((\tilde{H}_n - \tilde{H}_1)_{1 \le n < T_2} \in \cdot),$$

where $\tilde{H}$ [constructed from $(\tilde{\varepsilon}_t)_{0 \le t < F}$] is an alternating walk with exponential$(\theta - \mu)$ falls and exponential$(\theta + \mu)$ rises killed when reaching height $-F$.

Finally, for the distribution of $(\varepsilon_F(n))_{0 \le n \le 2\Gamma_F - 1}$ we avoid tedious calculations by using the (necessary) time-change argument on the whole path $(H_n)_{0 \le n \le 2N}$, $H'_n = H_{2N-n} - H_{2N}$. This maps the minimum on the minimum, rises on falls and vice versa, and the same argument can be applied: the process of excursions above the minimum is a Poisson point process with the fall parameter as its intensity, now $\theta + \lambda$, indexed by the ladder heights, killed when the absolute minimum is reached which now happens a.s. unless $\lambda = 0$. Note that in the calculations of the analogues of $r$ and $\sum_{n \ge 1} q^n \nu_n$, the component $\nu'_\infty$ of infinite excursions appears in a straightforward way because these excursions contain marks a.s. and do not belong to the pre-minimum path. Specifically, we calculate

$$1 - \nu'_\infty([0, \infty)^{\mathbb{N}}) = \sum_{1 \le n < \infty} C_n \left( \frac{\theta + \lambda}{2\theta} \right)^{n-1} \left( \frac{\theta - \lambda}{2\theta} \right)^n = \frac{\theta - \lambda}{\theta + \lambda}$$

$$\implies \quad \nu'_\infty([0, \infty)^{\mathbb{N}}) = \frac{2\lambda}{\theta + \lambda}$$

and

$$r' = \sum_{1 \le n \le \infty} (1 - q^n) \nu'_n([0, \infty)^{2n})$$

$$= \sum_{1 \le n \le \infty} \left( 1 - \left( \frac{\theta^2 - \mu^2}{\theta^2 - \lambda^2} \right)^n \right) C_n \left( \frac{\theta + \lambda}{2\theta} \right)^{n-1} \left( \frac{\theta - \lambda}{2\theta} \right)^n$$

$$= \theta + \lambda - (\theta + \mu) \frac{\theta - \mu}{\theta + \mu} = \mu + \lambda.$$

Note that we only use the walk-tree correspondence (Proposition 2) for finite trees; this extension can, for example, be obtained by density calculations.

The analogous argument for the law of the pre-minimum path completes the second proof of Theorem 2, since the point processes of excursions determine and identify the distributions of the pre- and post-minimum walks.

For the proof of the proposition, just apply the Markov property in $N = N_1$ and conclude by induction. $\quad \square$

We can now derive the forest composition rule from Proposition 4.

SECOND PROOF OF THEOREM 1. It will be convenient to switch to continuous time, take the Harris path $H^c$, rather than the walk, of a binary$(\lambda, \theta)$ forest and place independent marks at geometric$(1 - q)$ spaced leaves (local maxima). Denote their times in the path by $T^c_k$, $k \ge 1$, $T^c_0 = 0$, and also define the times of intertwining minima

$$S^c_k = \operatorname*{arg\,min}_{s \in [T^c_k, T^c_{k+1}]} H^c(s), \qquad k \ge 0.$$

Clearly, Williams' decomposition and Proposition 4 hold for Harris walks and paths with the obvious modifications of path laws. We deduce that the walk restricted to marked (black) leaves and intertwining minima is an exponential$(\mu \pm \lambda)$ alternating walk $H^b$ corresponding to a binary$(\lambda, \mu)$ forest. Its falls are now given by $F_k = H^c(T_k^c) - H^c(S_k^c)$, $k \geq 0$, its rises by $R_k = H^c(T_k^c) - H^c(S_{k-1}^c)$, $k \geq 1$.

Define the alternating lower envelope of $H^c$ derived from $T_k^c$, $k \geq 0$:

$$A_t = \begin{cases} \min\{H^c(s) : T_k^c \leq s \leq t\}, & \text{if } t \in [T_k^c, S_k^c], \\ \min\{H^c(s) : t \leq s \leq T_{k+1}^c\}, & \text{if } t \in [S_k^c, T_{k+1}^c]. \end{cases}$$

By Proposition 4, the excursions of $Z = H^c - A$ are (path excursions corresponding to the walk) excursions $(\bar{\varepsilon}_t)_{t \geq 0}$ of an alternating exponential$(\theta \pm \mu)$ Harris walk, separated by exponential$(\theta - \mu)$ spacings, $m \geq 1$. Put $T_0^b = 0$ and define inductively the successive times (in local time) of these points

$$T_n^b = \inf\left\{ t \geq T_{n-1}^b : Z\left( t + \sum_{m=1}^{n-1} L_m \right) > 0 \right\},$$

lengths of excursions

$$L_n = \inf\left\{ s > 0 : Z\left( T_n^b + \sum_{m=1}^{n-1} L_m + s \right) = 0 \right\},$$

and excursions

$$H^{(n)}(s) = H^c\left( T_n^b + \sum_{m=1}^{n-1} L_m + s \right) - H^c\left( T_n^b + \sum_{m=1}^{n-1} L_m \right), \qquad 0 \leq s \leq L_n.$$

We see that we have derived from $H^c$ the structure of Proposition 3 in terms of $H^b$, locations $(T_m^b)_{m \geq 0}$ and excursions $(H^{(m)})_{m \geq 1}$. We only need to check the independence of the three parts. Define $N_t = \sum_{n \geq 1} \mathbb{1}_{\{T_n^b \leq t\}}$ the arrival process of excursions (in local time). By Proposition 4 $(N_t)_{0 \leq t < F_0}$ is the arrival process of $\hat{\varepsilon}^{(0)}$, and is therefore a Poisson arrival process killed at the exponential fall $F_0$, and the same is true on each fall and rise of $H^b$. By the lack of memory property, $(N_t)_{0 \leq t < F_0 + R_1}$ is a Poisson arrival process killed at $F_0 + R_1$, and so on, and the independent killing time increases to infinity, and a limit argument establishes the independence of $(N_t)_{0 \leq t < \infty}$ and $H^b$. The same argument works for the Poisson point process of excursions with arrivals according to $N$. Note that the excursions attached to falls and rises have the same law by Lemma 4.

This proves not only Theorem 1 but also the stronger Proposition 1 since the coloring of geometric$(1 - q)$ spaced leaves means coloring every leaf (black) independently with probability $1 - q$.  $\square$

## 5. Williams' decomposition for Brownian motion with drift.

5.1. *From Harris paths to Brownian Williams decompositions.* Williams'
path decomposition for Brownian motion and related processes [28] is a
result which has received considerable attention by many authors [9, 11,
12, 16, 18, 20, 25]. In particular, Le Gall [18] uses discrete approximation
and tree arguments as well, but works with lattice approximations of three-
dimensional Bessel processes that arise by time reversal of Brownian motions
stopped at first passage times.

THEOREM 3 (Williams' decomposition [28]). *Let $(B_t)_{t \geq 0}$ be Brownian
motion with drift* $-\lambda$, $\lambda \in [0, \infty)$, *and* $T \sim$ *exponential$(\frac{1}{2}\kappa^2)$ independent.
Define the* (*absolute*) *height* $F = -\min\{B_t : 0 \leq t \leq T\}$ *and time* $M$, $F =
-B_M$, *of the minimum before* $T$, *and the terminal height* $R = B_T + F$ *above
the minimum.*

*Then the path* $(B_t)_{0 \leq t \leq T}$ *decomposes into two independent parts* $(B_t)_{0 \leq t \leq M}$
*and* $(B_{M+t} - B_M)_{0 \leq t \leq T-M}$. *Furthermore, if we put* $\mu^2 = \lambda^2 + \kappa^2$,

$$F \sim \text{exponential}(\mu - \lambda),$$

$$R \sim \text{exponential}(\mu + \lambda),$$

$$(B_t)_{0 \leq t \leq M} \sim (\tilde{B}_t)_{0 \leq t \leq \tau(-F)},$$

$$(B_{M+t} - B_M)_{0 \leq t \leq T-M} \sim (\tilde{B}_{\tau(-R)-t} - \tilde{B}_{\tau(-R)})_{0 \leq t \leq \tau(-R)},$$

*where* $\tilde{B}$ *is a Brownian motion with drift* $-\mu$ *independent of* $F$ *and* $R$, *and*
$\tau(-x) = \inf\{t \geq 0 : \tilde{B}_t < -x\}$.

PROOF. We take the limit $\theta \to \infty$ in Theorem 2: denote by $S_n^\theta = H_{2n}$
the random walk embedded in the Harris walk, its increments by $X_j^\theta =
-Y_j^\theta + Z_j^\theta$, $Y_j^\theta \sim$ exponential$(\theta - \lambda)$, $Z_j^\theta \sim$ exponential$(\theta + \lambda)$. Then we have

$$
\begin{aligned}
(1) \qquad \mathbb{P}\Big(\max_{1 \leq j \leq (1/2)\theta^2 n} |X_j^\theta| > \varepsilon\Big) &\leq 2(1 - (1 - e^{-\varepsilon(\theta-\lambda)})^{(1/2)\theta^2 n}) \\
&\leq \theta^2 n e^{-\varepsilon(\theta-\lambda)} \to 0.
\end{aligned}
$$

Note further that

$$\mathbb{E}(S_{[(1/2)\theta^2 t]}) = -\frac{2\lambda}{\theta^2 - \lambda^2}\Big[\frac{1}{2}\theta^2 t\Big] \quad \text{and} \quad \text{Var}(S_{(1/2)\theta^2 n}) = \frac{2(\theta^2 + \lambda^2)}{(\theta^2 - \lambda^2)^2}\Big[\frac{1}{2}\theta^2 t\Big].$$

Hence, Donsker's theorem (cf., e.g., Theorem 7.7.3 in [7]), implies that the
continuous-time process $S_{[(1/2)\theta^2 t]}^\theta$ converges weakly to Brownian motion
with drift $-\lambda$, as $\theta \to \infty$, uniformly on compact time intervals.

The same argument works if we replace $H_{2n}$ by $H_{2n+1} < H_{2n}$ since the in-
crement distribution is the same, and the negative starting value vanishes in
the limit. Furthermore, both processes converge simultaneously to the same

Brownian motion since (1) shows that their difference becomes uniformly small. Therefore, the intermediate process $H_{[\theta^2 t]}^\theta$ converges to Brownian motion with drift $-\lambda$, uniformly on compact time intervals.

Now, look at the decomposition of the Harris walks in the limit. The independent geometric times, after time scaling, converge as follows:

$$\mathbb{P}\left(\frac{2N_\theta}{\theta^2} > t\right) = \mathbb{P}\left(N_\theta > \left[\frac{1}{2}\theta^2 t\right]\right) = \left(1 - \frac{\mu^2 - \lambda^2}{\theta^2 - \lambda^2}\right)^{[(1/2)\theta^2 t]} \to e^{-(\mu^2 - \lambda^2)t/2}.$$

Hence, the distributional limit is exponential($\frac{1}{2}\kappa^2$) as required.

Now the decomposition can be read by taking limits in the following representation of Theorem 2:

$$\mathbb{E}(f(F^\theta)g((H_{[\theta^2 t]}^\theta)_{0 \leq t \leq M^\theta/\theta^2})h(R^\theta)k((H_{M^\theta+[\theta^2 t]}^\theta - H_{M^\theta}^\theta)_{0 \leq t \leq (2N^\theta - M^\theta)/\theta^2}))$$

$$= \mathbb{E}(f(F^\theta)g((H_{[\theta^2 t]}^\theta)_{0 \leq t \leq M^\theta/\theta^2}))$$

$$\times \mathbb{E}(h(R^\theta)k((H_{M^\theta+[\theta^2 t]}^\theta - H_{M^\theta}^\theta)_{0 \leq t \leq (2N^\theta - M^\theta)/\theta^2}))$$

for all bounded continuous functions $f, h$ and path functionals $g, k$. Specifically, $g = k = 1$ yields the exponential laws of $F$ and $R$, and taking $f = h = 1$ yields the distributions of the subpaths, employing the above argument and Donsker's theorem with different parameters. $\square$

5.2. *The forest growth process derived from Brownian motion.* Let $N$ be a homogeneous Poisson point process on $(0, \infty)^2$, with unit rate per unit area, assumed independent of a standard Brownian motion $B$. Consider the forest-valued process $(\mathcal{F}^\theta)_{\theta \geq 0}$ generated by sampling $B$ at times

$$T_n^\theta := \inf\{t \geq 0 : N([0, t] \times [0, \tfrac{1}{2}\theta^2]) = n\}, \qquad n \geq 0.$$

To be more precise, $\mathcal{F}^0$ is the trivial forest with no trees, just a forest floor identified with $[0, \infty)$, while for each $\theta > 0$ the random forest $\mathcal{F}^\theta$ is that associated with the Harris walk

$$H_{2n}^\theta := B_{T_n^\theta}, \qquad H_{2n+1}^\theta := \inf\{B_s : T_n^\theta < s < T_{n+1}^\theta\}, \qquad n \geq 0.$$

We call $(\mathcal{F}^\theta)_{\theta \geq 0}$ the *Poisson-sampled Brownian forest process*. We will exploit the well-known fact that this alternating walk has independent exponential step distributions. The following lemma is an expression in our forest framework of a characterization of Poisson-sampled Brownian trees due to Shapiro [27].

LEMMA 6.  *For each fixed $\theta > 0$, the forest $\mathcal{F}^\theta$ is a binary$(0, \theta)$ forest.*

PROOF. In view of Proposition 2, the coding of forests by alternating walks, this follows from the strong Markov property of $B$ at the times $T_n^\theta$, and Williams' decomposition of $B$ at the time $M_1^\theta$ of its minimum before the exponential($\frac{1}{2}\theta^2$) time $T_1^\theta$ (Theorem 3), whereby the random variables $-B_{M_1^\theta}$ and $B_{T_1^\theta} - B_{M_1^\theta}$ are independent exponential($\theta$) variables. See also [14, 25, 27] for variations of this argument. $\square$

The development of $(\mathcal{F}^\theta)_{\theta \geq 0}$ as $\theta$ varies is described by the next lemma, which is deduced using standard properties of Poisson processes:

LEMMA 7. *For each $0 < \lambda < \mu$, conditionally given $\mathcal{F}^\mu$, the forest $\mathcal{F}^\lambda$ is derived from $\mathcal{F}^\mu$ by taking the subforest of $\mathcal{F}^\mu$ spanned by a random set of leaves of $\mathcal{F}^\mu$ picked by a process of independent Bernoulli trials, where for each $i$ the $i$th leaf of $\mathcal{F}^\mu$ in order of depth-first search is put in the spanning set with probability $\lambda^2/\mu^2$.*

This lemma and Lemma 6 determine the joint distribution of $\mathcal{F}^\lambda$ and $\mathcal{F}^\mu$ for arbitrary $0 \leq \lambda < \mu$. This in turn determines the distribution of the whole forest-growth process $(\mathcal{F}^\theta)_{\theta \geq 0}$, because it turns out to be Markovian. This Markov property is not obvious, but a consequence of the following theorem.

THEOREM 4. *The Poisson-sampled Brownian forest process $(\mathcal{F}^\theta)_{\theta \geq 0}$ has independent growth increments, such that for each $0 \leq \lambda < \mu$ the forest of increments $\mathcal{F}^{\lambda,\mu}$ is a binary($\lambda, \mu$) forest.*

Note that Theorem 4 implies the forest composition rule (Theorem 1). We will prove Theorem 4 using Williams' decomposition (Theorem 3). Therefore we obtain a *third proof of Theorem 1*. This third proof is related to the second proof (using Williams' decomposition for Harris walks), but technically more involved, since the rigorous arguments in continuous time require more care.

From Theorem 4, one easily deduces the infinitesimal behavior.

COROLLARY 4. *The Poisson-sampled Brownian forest $(\mathcal{F}^\theta)_{\theta \geq 0}$ grows from the empty forest $\mathcal{F}^0$ by the following inhomogeneous infinitesimal transition rules:*

   (i) *at each time $\theta$, along each side of each edge of $\mathcal{F}^\theta$, twigs are attached to that side according to Poisson process with rate 1 per unit length of side per unit time;*

   (ii) *given that a twig is attached to a point $x$ on some side of $\mathcal{F}^\theta$, the length of that twig has an exponential($2\theta$) distribution.*

This twig-growth description of the Poisson-sampled Brownian forest is almost the same as Aldous' [2] process 2 for construction of a self-similar continuum random tree. Aldous' construction differs from ours by a rotation of 90° which turns our horizontal forest floor into an infinitely high spine of a tree, and a duplication which allows an independent forest with same distribution to grow on the right-hand side of the spine.

To prove Theorem 4, the following quantities, related to the time reversal in Williams' decomposition, will be useful.

Given a continuous path $B = (B_t)_{t \geq 0}$ and an increasing sequence of sampling times $(T_n)_{n \geq 1}$, as in the second proof of Theorem 1, we define the *alternating lower envelope of $B$ derived from* $(T_n)_{n \geq 1}$ to be the process with locally bounded variation $A_t := A_t(B, (T_n)_{n \geq 1})$ constructed as follows:

$$A_t = \begin{cases} \min\{B_s : T_n \leq s \leq t\}, & \text{if } t \in [T_n, M_{n+1}], \\ \min\{B_s : t \leq s \leq T_{n+1}\}, & \text{if } t \in [M_{n+1}, T_{n+1}], \end{cases}$$

where $T_0 = 0$ and $M_{n+1}$ is a time in $[T_n, T_{n+1}]$ at which $B$ attains its minimum on that interval. Note that $M_{n+1} = T_n$ or $M_{n+1} = T_{n+1}$ may occur here, causing *degeneracy* in the sequel. We therefore assume throughout that $B$ and $T_n$ are such that this does not happen, as is the case almost surely in our applications to Poisson-sampled Brownian paths.

Call the process $Z_t := B_t - A_t \geq 0$ the *reflected process* derived from $B$ and $(T_n)_{n \geq 1}$. Define a random sign process $\sigma_t$, with $\sigma_t = -1$ if $t$ is in one of the intervals $(T_n, M_{n+1}]$ when $A$ is decreasing, and $\sigma_t = +1$ if $t$ is in one of the intervals $(M_{n+1}, T_{n+1}]$ when $A$ is increasing, so the process

$$J_t := \int_0^t \sigma_s \, dA_s$$

is a continuous increasing process, call it the *increasing process derived from $B$ and* $(T_n)_{n \geq 1}$.

The following (deterministic) lemma now follows easily. It is a generalization of Proposition 3 to composite forests derived from continuous paths.

LEMMA 8. *For an arbitrary continuous path $B$ with $\inf\{B_t : t \geq 0\} = -\infty$, let $\mathcal{F}^b$ be the* black *forest derived from $B$ by sampling at some increasing sequence of (black) times $(T_n^b)_{n \geq 1}$ (not leading to degeneracy), and let $\mathcal{F}^c$ be the forest derived from $B$ by sampling at some increasing sequence of times $(T_n^c)_{n \geq 1}$, where $\{T_n^c : n \geq 1\} = \{T_n^b : n \geq 1\} \cup \{T_n^r : n \geq 1\}$ for some increasing sequence of (red) times $(T_n^r)_{n \geq 1}$ (again not leading to degeneracy). Let $A_t = A_t(B, (T_n^b)_{n \geq 1})$ be the alternating lower envelope of $B$ induced by the $(T_n^b)_{n \geq 1}$, let $Z := B - A$ be the reflected process and let $J_t = J_t(B, (T_n^b)_{n \geq 1})$ be the increasing process derived from $B$ and the $(T_n^b)_{n \geq 1}$. Then the forest $\mathcal{F}^r$ of red innovations grown onto $\mathcal{F}^b$ to form $\mathcal{F}^c$ is identical to the forest derived from $Z - J$ by sampling at the times $(T_n^r)_{n \geq 1}$.*

Theorem 4 now follows from the construction of the forest growth process $(\mathcal{F}^{\theta})_{\theta \geq 0}$ by Poisson sampling of $B$, the previous lemma and the following reformulation of Williams' decomposition.

LEMMA 9. *For $\lambda \geq 0$ let $\mathbb{P}_{-\lambda}$ govern $(B_t)_{t \geq 0}$ as a Brownian motion with drift $-\lambda$, meaning that the $\mathbb{P}_{-\lambda}$ distribution of $(B_t)_{t \geq 0}$ is the $\mathbb{P}_0$ distribution of $(B_t - \lambda t)_{t \geq 0}$. Let $\mathbb{P}_{-\lambda}$ also govern $T_m^{\kappa}$ for $m \geq 1$ as the points of a Poisson process with rate $\frac{1}{2}\kappa^2$ which is independent of $B$. For $m \geq 1$ let $F_m^{\kappa} := B(T_{m-1}^{\kappa}) - B(M_m^{\kappa})$ be the $m$th fall and $R_m^{\kappa} := B(T_m^{\kappa}) - B(M_m^{\kappa})$ the $m$th rise of the alternating walk defined by the values $B(T_m^{\kappa})$ and the intermediate minima $B(M_m^{\kappa})$. Then under $\mathbb{P}_{-\lambda}$ for each $-\lambda \leq 0$:*

(i) *the random variables $F_m^{\kappa}$ and $R_m^{\kappa}$ are independent, with*

$$F_m^{\kappa} \sim \text{exponential}(\sqrt{\kappa^2 + \lambda^2} - \lambda),$$
$$R_m^{\kappa} \sim \text{exponential}(\sqrt{\kappa^2 + \lambda^2} + \lambda).$$

(ii) *If $A^{\kappa}$ is the alternating lower envelope derived from $B$ and the sample times $(T_m^{\kappa})_{m \geq 1}$, $Z^{\kappa}$ is the reflected process and $J^{\kappa}$ the increasing process, then the process $Z^{\kappa} - J^{\kappa}$ is a Brownian motion with drift $-\sqrt{\kappa^2 + \lambda^2}$. Equivalently, $Z^{\kappa}$ is a Brownian motion with drift $-\sqrt{\kappa^2 + \lambda^2}$ on $(0, \infty)$ and simple reflection at $0$.*

(iii) *The Brownian motion $Z^{\kappa} - J^{\kappa}$ with drift $-\sqrt{\kappa^2 + \lambda^2}$ is independent of the bivariate sequence of falls and rises $(F_m^{\kappa}, R_m^{\kappa})_{m \geq 1}$, hence also independent of the Poisson-sampled Brownian forest $\mathcal{F}^{\kappa}$ which they encode.*

PROOF. The argument parallels the second proof of Theorem 1. The independence assertions in (i), and the exponential form of the distributions of the falls and rises follow from Williams' decomposition (Theorem 3) at the time $M_1^{\kappa}$ of its minimum on the interval $[0, T_1^{\kappa}]$, and repeated application of the strong Markov property of $B$ at the times $T_m^{\kappa}$. Also according to Theorem 3, conditionally given $F_1^{\kappa} = f$ and $R_1^{\kappa} = r$, the fragments of the path of $B$ on the intervals $[0, M_1^{\kappa}]$ and $[M_1^{\kappa}, T_1^{\kappa}]$ are independent, the first fragment distributed like a Brownian motion with drift $-\sqrt{\kappa^2 + \lambda^2}$, started at $0$ and run until it first hits $-f$, while the second fragment reversed is like a Brownian motion with the same negative drift started at $r - f$ and run until it first hits $-f$. It follows that with the same conditioning, the two reflected path fragments $(Z_t^{\kappa})_{0 \leq t \leq T_1^{\kappa}}$ and $(Z_{M_1^{\kappa} + u}^{\kappa})_{0 \leq u \leq T_1^{\kappa} - M_1^{\kappa}}$ are independent, the first fragment a reflected Brownian motion with drift $-\sqrt{\kappa^2 + \lambda^2}$, run until its local time at $0$ reaches $f$, and the second fragment reversed a reflected Brownian motion with drift $-\sqrt{\kappa^2 + \lambda^2}$, run until its local time at $0$ reaches $r$. But from this description, and the well-known reversibility of a one-dimensional diffusion stopped at an inverse local time, still conditioning on $F_1^\kappa = f$ and $R_1^\kappa = r$, the process $(Z_t^\kappa)_{0 \le t \le T_1^\kappa}$ is identified as a reflected Brownian motion with drift $-\sqrt{\kappa^2 + \lambda^2}$ run until its local time at 0 first reaches $f + r$. Now by repeated use of this argument, the entire process $(Z_t^\kappa)_{t \ge 0}$ conditional on all the rises and falls is a reflected Brownian motion with drift $-\sqrt{\kappa^2 + \lambda^2}$ run forever, independent of the given values of the rises and falls, provided they sum to $\infty$ which they obviously do almost surely. Since $Z^\kappa - J^\kappa$ is by construction the Brownian motion driving this reflected process, the conclusions (ii) and (iii) are evident. $\square$

**Acknowledgments.** M. Winkel would like to thank the Department of Statistics at U.C. Berkeley for hospitality. This work was started at the St. Flour Summer School in July 2002, and completed in March 2004 during visits of the first author to the Laboratoire de Probabilités at the University of Paris 6, and of the second author to the Laboratoire Mathématiques at the University of Paris-Sud. We also thank the anonymous referees for several suggestions leading to improvements of the presentation.

## REFERENCES

[1] ABRAHAM, R. (1992). Un arbre aléatoire infini associé à l'excursion Brownienne. *Séminaire de Probabilités XXVI. Lecture Notes in Math.* **1526** 374–397. Springer, Berlin. MR1232004

[2] ALDOUS, D. (1991). The continuum random tree. I. *Ann. Probab.* **19** 1–28. MR1085326

[3] ALDOUS, D. (1993). The continuum random tree. III. *Ann. Probab.* **21** 248–289. MR1207226

[4] BORODIN, A. N. and SALMINEN, P. (1996). *Handbook of Brownian Motion—Facts and Formulae.* Birkhäuser, Basel. MR1477407

[5] DUQUESNE, T. and LE GALL, J.-F. (2002). *Random Trees, Lévy Processes and Spatial Branching Processes.* Société Mathématique de France, Paris.

[6] DUQUESNE, T. and WINKEL, M. (2005). Growth of Lévy forests. Unpublished manuscript.

[7] DURRETT, R. (1996). *Probability: Theory and Examples.* Duxbury Press, Belmont, CA. MR1609153

[8] EVANS, S. N., PITMAN, J. and WINTER, A. (2005). Rayleigh processes, real trees, and root growth with re-grafting. *Probab. Theory Related Fields.* To appear.

[9] FITZSIMMONS, P. J. (1986). Another look at Williams' decomposition theorem. In *Seminar on Stochastic Processes, Progress in Probability and Statistics* **12** 79–85. Birkhäuser, Boston. MR896736

[10] GEIGER, J. and KERSTING, G. (1997). Depth-first search of random trees, and Poisson point processes. In *Classical and Modern Branching Processes* (K. B. Athreya and P. Jagers, eds.) 111–126. Springer, New York. MR1601713

[11] GREENWOOD, P. and PITMAN, J. (1980). Fluctuation identities for random walk by path decomposition at the maximum. *Adv. in Appl. Probab.* **12** 291–293. MR588409

[12] GREENWOOD, P. and PITMAN, J. (1980). Fluctuation identities for Lévy processes and splitting at the maximum. *Adv. in Appl. Probab.* **12** 893–902. MR588409

[13] HARRIS, T. E. (1952). First passage and recurrence distributions. *Trans. Amer. Math. Soc.* **73** 471–486. MR52057

[14] Hobson, D. G. (2000). Marked excursions and random trees. *Séminaire de Probabilités XXXIV. Lecture Notes in Math.* **1729** 289–301. Springer, Berlin. MR1768069

[15] Kendall, D. G. (1951). Some problems in the theory of queues (with discussion). *J. Roy. Statist. Soc. Ser. B* **13** 151–185. MR47944

[16] Kersting, G. and Memisoglu, K. (2004). Path decompositions for Markov chains. *Ann. Probab.* **32** 1370–1390. MR2060301

[17] Knuth, D. E. (1969). *The Art of Computer Programming* **1**. *Fundamental Algorithms.* Addison–Wesley, Reading, MA. MR378456

[18] Le Gall, J.-F. (1986). Une approche élémentaire des théorèmes de décomposition de Williams. *Séminaire de Probabilités XX. Lecture Notes in Math.* **1204** 447–464. Springer, Berlin.

[19] Le Gall, J.-F. (1989). Marches aléatoires, mouvement brownien et processus de branchement. *Séminaire de Probabilités XXIII. Lecture Notes in Math.* **1372** 258–274. Springer, Berlin. MR1022916

[20] Millar, P. W. (1978). A path decomposition for Markov processes. *Ann. Probab.* **6** 345–348. MR461678

[21] Neveum, J. (1986). Arbres et processus de Galton–Watson. *Ann. Inst. H. Poincaré Probab. Statist.* **22** 199–207. MR850756

[22] Neveu, J. (1986). Erasing a branching tree. In *Analytic and Geometric Stochastics*: *Papers in Honour of G. E. H. Reuter* (*Special supplement to Adv. in Appl. Probab.*) 101–108. MR868511

[23] Neveu, J. and Pitman, J. (1989). Renewal property of the extrema and tree property of the excursion of a one-dimensional Brownian motion. *Séminaire de Probabilités XXIII. Lecture Notes in Math.* **1372** 239–247. Springer, Berlin. MR1022914

[24] Neveu, J. and Pitman, J. W. (1989). The branching process in a Brownian excursion. *Séminaire de Probabilités XXIII. Lecture Notes in Math.* **1372** 248–257. Springer, Berlin.

[25] Pitman, J. (1999). Brownian motion, bridge, excursion, and meander characterized by sampling at independent uniform times. *Electron. J. Probab.* **4** 1–33. MR1690315

[26] Pitman, J. (2002). Combinatorial stochastic processes. Lecture Notes for St. Flour Course July 2002. Technical Report 621, Berkeley.

[27] Shapiro, J. W. (1995). Capacity of Brownian trace and level sets. Ph.D. dissertation, Univ. California, Berkeley.

[28] Williams, D. (1974). Path decomposition and continuity of local time for one-dimensional diffusions. I. *Proc. London Math. Soc. (3)* **28** 738–768. MR350881

Department of Statistics #3860
University of California, Berkeley
367 Evans Hall
Berkeley, California 94720–3860
USA
e-mail: pitman@stat.berkeley.edu

Department of Statistics
University of Oxford
1 South Parks Road
Oxford OX1 3TG
United Kingdom
e-mail: winkel@stats.ox.ac.uk
url: www.stats.ox.ac.uk/˜winkel/indexe.html